\documentclass[11pt]{article}
\usepackage{amssymb}
\usepackage{latexsym}
\usepackage{epsfig}
\usepackage{psfrag}
\def\scfig #1 #2 {\resizebox{#2}{!}{\includegraphics{#1}}}

\newtheorem{Theorem}{\bf Theorem}
\newtheorem{Lemma}[Theorem]{\bf Lemma}
\newtheorem{Proposition}[Theorem]{\bf Proposition}
\newtheorem{Corollary}[Theorem]{\bf Corollary}

\newtheorem{Remark}[Theorem]{\bf Remark}

\def\qed{\hfill$\Box$}

\begin{document}

\title{A natural representation model for symmetric groups}
\author{{\sc Vijay Kodiyalam$^1$ and D.-N. Verma$^2$}\\
I.M.Sc.,
Chennai$^1$ and T.I.F.R., Mumbai$^2$\\
INDIA}

\maketitle

\begin{abstract}
For any positive integer $N$, we describe a natural complex representation of 
the symmetric group $\Sigma_N$ on the vector space spanned by its involutions
that contains each irreducible representation exactly once.
\end{abstract}

\section{Introduction}

Given a finite group $G$, it is interesting to construct a
`natural representation' of $G$ that contains each of its irreducible 
representations exactly once.
Such a representation is called a Gelfand model for $G$ in \cite{Stn} and 
constructed for the general linear groups over finite fields in \cite{Kly}.
In \cite{GdoRjo}, such a model is described for the family of symmetric groups
where the representation space is obtained as the zeroes of certain 
differential operators in the Weyl algebra.
In this note, we contruct another very simple model for this family.
It would be interesting to explicitly compare the two models and give a 
construction of either model for, say, the finite Weyl groups.

Throughout this paper we adopt the following notations and conventions.
Let $N$ be a fixed positive integer and $G = \Sigma_N$ be the symmetric group
of permutations of the set $[N] := \{1,2,\cdots,N\}$.
The identity element of $G$ will be denoted $1_G$.
For $0 \leq j \leq \lfloor N/2 \rfloor$, let $X_j$ be the conjugacy class of
involutions - elements with square $1_G$ - of $G$ that are a product of $j$ 
disjoint transpositions.
If $\tau \in X_j$ we will speak of $\tau$ as having length $j$ and write 
$\ell(\tau) = j$.
Let $X = \coprod _{j = 0}^{\lfloor N/2 \rfloor} X_j$ be the set of all 
involutions of $G$.
All our vector spaces and representations will be defined over the field of 
complex 
numbers.
Let $V_j$ be the vector space with basis $X_j$.

Each $V_j$ admits a representation $\pi_j$ of $G$ defined on the basis $X_j$ 
as follows: take
$\tau \in X_j$ and write it as $(p_1 p_2)(p_3 p_4) \cdots 
(p_{2j-1} p_{2j})$ where $p_{2k-1} < p_{2k}$ for each $k = 1,2,\cdots,j$.
For $\sigma \in G$, define:
$$
\pi_j(\sigma)(\tau) = S(\sigma,\tau)~(\sigma({p_1}) 
\sigma({p_2}))(\sigma({p_3}) \sigma({p_4})) \cdots 
(\sigma({p_{2j-1}}) \sigma({p_{2j}}))
$$
where $S(\sigma,\tau) \in \{\pm 1\}$ is defined to be ${(-1)}^{\#\{k  : 
\sigma(p_{2k-1}) > \sigma(p_{2k})\}}$ - where {\small{\#}} stands for 
`cardinality of'.
Equivalently, in matrix terms, we have that $\pi_j(\sigma)_\tau^{\tau^\prime} 
= 
S(\sigma,\tau) \delta_{\tau^\prime,\sigma\tau\sigma^{-1}}$.

In order to see that $\pi_j$ is indeed a representation of $G$, we note that
an alternative description is as follows.
Let $V$ denote the `natural' representation space of $\Sigma_N$ with a basis
$v_1,\cdots,v_N$ on which $\Sigma_N$ acts as $\sigma(v_k) = v_{\sigma(k)}$.
Consider the representation $\bigwedge^2 V$ which has basis $v_{ij} := v_i 
\wedge v_j$ for $1 \leq i < j \leq N$.
Identifying $Sym(\bigwedge^2 V)$ with the polynomial ring in the variables 
$v_{ij}$, it is easy to see that $V_j$ is isomorphic to the $j^{th}$ graded 
component of 
the $G$-module $Sym(\bigwedge^2 V)/I$ where $I$ is the monomial ideal generated
by all $v_{ij}v_{kl}$ where $\{i,j\} \cap \{k,l\} \neq \emptyset$.
This description has the virtue of providing a $G$-module algebra structure on
$A = \oplus_{j = 1}^{\lfloor N/2 \rfloor} V_j$.

Observe that the fact that $\pi_j$ is a representation of $G$ implies that
the function $S$ satisfies the following basic equation:
\begin{equation}\label{basic}
S(\sigma^\prime \sigma, \tau) = S(\sigma^\prime, 
\sigma\tau\sigma^{-1})S(\sigma, \tau),
\end{equation}
for any $\sigma,\sigma^\prime \in G$ and $\tau \in X$ and is, in particular,
a character of the centraliser of $\tau$.

We may now state our main theorem.

\begin{Theorem}\label{main}
The $G$-module $A = \oplus_{j = 1}^{\lfloor N/2 \rfloor} V_j$ contains every
irreducible representation of $G$ exactly once.
\end{Theorem}

The strategy of proof is very simple. We first show that the representations
$\pi_j$ are mutually disjoint in the sense that they do not contain any
non-zero isomorphic $G$-subrepresentations.
Then we show that the $G$-endomorpism ring of any $V_j$ is abelian thus
establishing that it is multiplicity free.
Finally we show that the dimension of $End_G(A)$ is equal to the number of
partitions of $N$ showing that every irreducible representation of $G$
features in $A$.

%Finally, we give a direct proof of the fact that all irreducible 
%representations of $G$ are defined over ${\mathbb Q}$ and then appeal to a 
%classical
%count of involutions due to Frobenius and Schur to see that their number
%agrees with the sum of the dimensions of its
%irreducible representations.

Much of the calculations are based on analysing matrix entries of a $G$-map
from $V_j$ to $V_k$ with respect to the distinguished bases $X_j$ and $X_k$.
So let $T : V_j \rightarrow V_k$ be a $G$-map. We will identify $T$ with
its matrix $((T_\tau^\kappa))$ for $\tau \in X_j$ and $\kappa \in X_k$.
Thus $T(\tau) = \Sigma_{\tau} T_\tau^\kappa \kappa$. That $T$ is a $G$-map
is seen to translate into the equations:
\begin{eqnarray}\label{teqs}
S(\sigma,\kappa) T^\kappa_\tau = S(\sigma,\tau) 
T^{\sigma\kappa\sigma^{-1}}_{\sigma\tau\sigma^{-1}}
\end{eqnarray}
for all $\tau \in X_j$ and $\kappa \in X_k$.

\section{Involutions in the symmetric group}

To proceed, we will need to analyse the structure of the set $X \times X$ as a
$G$-set under the simultaneous conjugation action.
Given an element $(\tau,\kappa) \in X \times X$, we define its characteristic
partition, denoted $P(\tau,\kappa)$,  to be the finest partition, say
$\{H_1,H_2,\cdots,H_t\}$, of $[N]$
such that both $\tau$ and $\kappa$ are in the subgroup $\Sigma_{H_1} \times
\cdots \times \Sigma_{H_t}$ (with the obvious meaning)  of $\Sigma_N$.
There are unique expressions $\tau = \tau_1\tau_2\cdots\tau_t$ and $\kappa =
\kappa_1\kappa_2\cdots\kappa_t$ where $\tau_i,\kappa_i \in \Sigma_{H_i}$
for $i = 1,2,\cdots,t$.
This also specifies a partition $p(\tau,\kappa)$ of the integer $N$ as
$\Sigma_{i=1}^{t} |H_i|$ which we will refer to as the numerical
characteristic partition of $(\tau,\kappa)$.
Clearly $p(\tau,\kappa)$ is an invariant of the orbit of $(\tau,\kappa)$.

\begin{Lemma}\label{lemma2} Let $(\tau,\kappa) \in X \times X$ be such that 
$P(\tau,\kappa) = [N]$. Then,\\
(a) The following conditions are equivalent:\\
%\begin{itemize}
{~~(i)} $\ell(\tau) = \ell(\kappa)$, and\\
{~~(ii)} There is an involution $\lambda \in \Sigma_N$ that conjugates each of 
$\tau$ and
$\kappa$ to the other.\\
%\end{itemize}
(b) If the conditions (i) and (ii) hold, then $N - 2\ell(\tau) \leq 1$. For
any $\sigma \in \Sigma_N$ that commutes with both $\tau$ and $\kappa$,
$S(\sigma,\tau)S(\sigma,\kappa) = 1$.\\
(c) If the conditions (i) and (ii) do not hold, there is an involution $\sigma 
\in \Sigma_N$ commuting with both $\tau$ and $\kappa$ and is such that 
$S(\sigma,\tau)S(\sigma,\kappa) = -1$.
\end{Lemma}

\noindent
{\bf Proof.}
Suppose that there is an element $p_1$ of $[N]$ that is moved by one
of the involutions, say $\tau$, to $p_2$ and fixed by $\kappa$.
Then, unless $N=2$, $p_2$ is neccessarily moved by $\kappa$ to $p_3$ for
otherwise $H = \{p_1,p_2\}$ would be a non-empty subset of $[N]$ with 
non-empty complement $K$ so that
both $\tau$ and $\kappa$ are in $\Sigma_H \times \Sigma_K$.
Then, unless $N=3$, $p_3$ is neccessarily moved by $\tau$, to  $p_4$, for
otherwise, $H = \{p_1,p_2,p_3\}$ would be a non-empty subset of $[N]$ with 
non-empty complement $K$ so that
both $\tau$ and $\kappa$ are in $\Sigma_H \times \Sigma_K$.
Proceeding in this manner, we see that one of the following must hold:\\
Case I: $N$ is even and there is a renumbering $\{p_1,p_2,\cdots,p_{2j}\}$ of 
$[N]$ with $\tau = (p_1 
p_2)(p_3 p_4) \cdots 
(p_{2j-1} p_{2j})$ and $\kappa = (p_2 p_3)(p_4 p_5) \cdots 
(p_{2j-2} p_{2j-1})$, or\\
Case II: 
$N$ is odd and there is a renumbering  $\{p_1,p_2,\cdots,p_{2j+1}\}$ of 
$[N]$ with $\tau = (p_1 
p_2)(p_3 p_4) \cdots 
(p_{2j-1} p_{2j})$ and $\kappa = (p_2 p_3)(p_4 p_5) \cdots 
(p_{2j} p_{2j+1})$.

On the other hand, if no element of $[N]$ is left fixed by either $\kappa$
or $\tau$,
then $N$ is clearly neccessarily even and
a similar analysis shows that  the following holds:\\
Case III: $N$ is even and there is a renumbering 
$\{p_1,p_2,\cdots,p_{2j}\}$ of $[N]$ with $\tau = 
(p_1 
p_2)(p_3 p_4) \cdots 
(p_{2j-1} p_{2j})$
and $\kappa = (p_2 p_3)(p_4 p_5) \cdots 
(p_{2j-2} p_{2j-1})(p_{2j} p_1)$.

To prove (a), since (ii) clearly implies (i), now suppose that (i) holds. Then 
we are in either Case II or Case 
III.
In Case II, the involution $\lambda$ that interchanges $p_i$ with 
$p_{N+1-i}$ conjugates each of
$\tau$ and
$\kappa$ to the other while in Case III, the involution $\lambda$ that 
interchanges $p_i$ with $p_{N+2-i}$ works (the addition is mod N).

The first assertion of (b) is clear since in Case II, $N = 2j +1$ while in 
Case III, $N 
= 2j$ with $j = \ell(\tau) = \ell(\kappa)$. The second assertion is a little 
more delicate. Observe that as $P(\tau,\kappa) = [N]$, any $\sigma$ commuting 
with both $\tau$ and
$\kappa$ is determined by any one $\sigma(k)$. In Case II, since $p_1$ is
the unique element of $[N]$ moved by $\tau$ and fixed by $\kappa$, $\sigma(p_1)
= p_1$ and hence $\sigma = 1_G$.

In Case III, the elements that commute with both $\tau$ and $\kappa$ form a 
dihedral group generated by the involution
$\sigma_2 = (p_1~p_{2j})(p_2~p_{2j-1})~\cdots(p_j~p_{j+1})$ and the element 
$\sigma_1 = 
(p_1~p_3~\cdots~p_{2j-1})(p_2~p_4~\cdots~p_{2j})$. Since $S$ is a character on 
the centraliser
of $\tau$, it suffices to check that  $S(\sigma_i,\tau)S(\sigma_i,\kappa) = 1$.
Inspection shows that $S(\sigma_1,\kappa) = 1 = S(\sigma_1,\tau)$ and that
$S(\sigma_2,\kappa) = (-1)^j = S(\sigma_2,\kappa)$, finishing the proof of (b).

The hypothesis of (c) implies that we are in Case I and here again the 
involution that interchanges $p_i$ with $p_{N+1-i}$ is seen to satisfy the 
conclusion of (c). \qed

\begin{Corollary} \label{corr} Let $(\tau_1,\kappa_1),(\tau_2,\kappa_2) \in X 
\times X$ 
with $P(\tau_1,\kappa_1) = [N] = P(\tau_2,\kappa_2)$ and with $\ell(\tau_1)
= \ell(\kappa_1)$ and  $\ell(\tau_2)
= \ell(\kappa_2)$.
Then there is a $\sigma \in \Sigma_N$ such that $(\tau_2,\kappa_2) = 
(\sigma \tau_1\sigma^{-1},\sigma \kappa_1 \sigma^{-1})$.
\end{Corollary}

\noindent
{\bf Proof.} It follows from Lemma \ref{lemma2} that for each of 
$(\tau_1,\kappa_1)$ and $(\tau_2,\kappa_2)$ we are either Case II or
Case III, and in the same case for both since that depends only on whether
$N$ is odd or even. Thus from the proof of that lemma, $(\tau_1,\kappa_1)$ and 
$(\tau_2,\kappa_2)$ differ only by a renumbering of $[N]$ which yields the
sought after $\sigma$.\qed

\begin{Lemma}\label{invlemma}
Let $(\tau,\kappa) \in X \times X$ with
$P(\tau,\kappa) = \{H_1,\cdots,H_t\}$.
Then exactly one of the following two conditions holds:\\
(i) There ia a  $\sigma \in \Sigma_N$  commuting with $\tau$ 
and 
$\kappa$ such that
$S(\sigma,\tau)S(\sigma,\kappa) = -1$.\\
(ii) There is an involution $\lambda \in \Sigma_{H_1} \times \cdots \times 
\Sigma_{H_t}$ such that $\kappa =
\lambda\tau\lambda^{-1}$.
\end{Lemma}

\noindent
{\bf Proof.} To see that at least one of the two conditions holds, consider the
elements $\tau_i$ and $\kappa_i$. If their lengths are equal for all $i = 
1,\cdots,t$, 
then by Lemma \ref{lemma2}, they are conjugate by involutions $\lambda_i \in 
\Sigma_{H_i}$. The product $\lambda = \lambda_1 \cdots \lambda_t$ gives 
an involution
living in $\Sigma_{H_1} \times \cdots \times 
\Sigma_{H_t}$ that conjugates $\tau$ to $\kappa$. Otherwise,
for some $i$ there is, by Lemma \ref{lemma2} again, a $\sigma_i \in 
\Sigma_{H_i}$ that commutes with both 
$\tau_i$
and $\kappa_i$ and is such that, in $\Sigma_{H_i}$, 
$S(\sigma_i,\tau_i)S(\sigma_i,\kappa_i) = -1$.
Simply regard $\sigma = \sigma_i$ as an element of $\Sigma_N$ to get an element
satisfying (i).

To see that they are mutually exclusive,
suppose that (ii) holds and consider any $\sigma \in \Sigma_N$ commuting with 
both
$\tau$ and $\kappa$.
Then $\sigma$ preserves the characteristic partition in the sense that 
$\{\sigma(H_1),\cdots,\sigma(H_t)\} = \{H_1,\cdots,H_t\}$.
Recalling the definition of $S$, a little thought shows that $S(\sigma,\tau)$
is the product of all $S(\sigma_i,\tau_i)$
where the product is over all $i$ for which $\sigma(H_i) = H_i$ and $\sigma_i 
= \sigma|_{H_i}$.
A similar statement holds for $S(\sigma,\kappa)$. We now note that $\sigma_i$
commutes with $\tau_i$ and $\kappa_i$ and appeal to Lemma \ref{lemma2}(b)
to finish the proof.\qed

We note that a consequence of this lemma is that for any $(\tau,\kappa) \in X 
\times X$, whether the condition (i) or (ii) of the lemma holds is a property
of its $G$-orbit - a fact that does not seem to be directly verifiable for the
condition (i).

\section{The main theorem}

\begin{Proposition}\label{disjoint}
For $j \neq k$ the representations $\pi_j$ and $\pi_k$ are mutually disjoint.
\end{Proposition}

\noindent
{\bf Proof.}
It is to be seen that any $G$-map $T: V_j \rightarrow V_k$ is zero. Since the
entries of $T$ satisfy the system of equations (\ref{teqs}), it suffices to
see that for any $\tau \in X_j$ and $\kappa \in X_k$ there exists a $\sigma
\in G$ commuting with both so that $S(\sigma,\tau)$ and $S(\sigma,\kappa)$
have different signs. But that follows from Lemma \ref{invlemma} since the
case (ii) there cannot occur for involutions of different length.\qed

\begin{Proposition}\label{comm}
For any $j$, the endomorphism ring $End_G(V_j)$ is commutative.
\end{Proposition}

\noindent
{\bf Proof.}
Let $T$ be an arbitrary element of $End_G(V_j)$. We assert that the
matrix of $T$ with respect to the basis $X_j$ is symmetric.
To see this, let $\tau,\tau^\prime$ be arbitrary elements of $X_j$.
If there is an element $\sigma \in \Sigma_N$ commuting with both so that
$S(\sigma,\tau)S(\sigma,\tau^\prime) = -1$, then equations (\ref{teqs})
imply that both $T^\tau_{\tau^\prime}$ and $T_\tau^{\tau^\prime}$ are $0$.
Otherwise, Lemma \ref{invlemma} implies that there is an involution $\lambda 
\in \Sigma_N$
that conjugates each of $\tau$ and $\tau^\prime$
to the other. By the basic property of the $S$ function, $1 = S(1_G,\tau)
= S(\lambda,\tau)S(\lambda^{-1},\lambda\tau\lambda^{-1}) = 
S(\lambda,\tau)S(\lambda,\kappa)$. Now again appeal to the equations 
(\ref{teqs})
to see that $T^\tau_{\tau^\prime} = T_\tau^{\tau^\prime}$ in this case too.

But now, $End_G(V_j)$ is identified with an algebra of symmetric matrices
and is therefore commutative.\qed

\begin{Proposition}\label{dimn}
For any $j$, the dimension of $End_G(V_j)$ is equal to the number of partitions
of $N$ with $N-2j$ odd parts.
\end{Proposition}

\noindent
{\bf Proof.} It should be clear from the equations (\ref{teqs}) that the matrix
entries of an arbitrary $T \in End_G(V_j)$ are completely determined by the
values of $T^{\tau^\prime}_{\tau}$ as $(\tau,\tau^\prime)$ range over 
representatives of  simultaneous conjugation $G$-orbits in $X_j \times X_j$. 

Furthermore, if for a $(\tau,\tau^\prime) \in X_j \times X_j$, there exists 
$\sigma \in G$ commuting with both such that 
$S(\sigma,\tau)S(\sigma,\tau^\prime) = -1$ then the entries of $T$ in the
orbit of $(\tau,\tau^\prime)$ vanish as observed earlier.
Conversely, suppose that it is a consequence of equations (\ref{teqs}) that
$T^{\tau^\prime}_{\tau} = 0$ for all $T \in End_G(V_j)$. Then neccessarily 
there
exist $(\tau_i,\tau^\prime_i)$ for $i = 1,2,\cdots,n$ in the orbit of 
$(\tau,\tau^\prime)$ and $\sigma_i \in G$ so that for $i = 1,2,\cdots,n-1$, we 
have
$(\tau_{i+1},\tau_{i+1}^\prime) = \sigma_i\cdot (\tau_{i},\tau_{i}^\prime)$
and $(\tau_{1},\tau_{1}^\prime) = \sigma_n\cdot (\tau_{n},\tau_{n}^\prime)$
and such that 
$$
\frac{S(\sigma_n,\tau_n^\prime)S(\sigma_{n-1},\tau_{n-1}^\prime)\cdots 
S(\sigma_1,\tau_1^\prime)}
{S(\sigma_n,\tau_n)S(\sigma_{n-1},\tau_{n-1})\cdots S(\sigma_1,\tau_1)} = -1.
$$
But now the product $\sigma = \sigma_n\cdots \sigma_1$ is easily verified - 
using equation (\ref{basic}) - to commute with both and satisfy 
$S(\sigma,\tau_1)S(\sigma,\tau_1^\prime) = -1$ .
Therefore the 
dimension
of $End_G(V_j)$ is equal to the number of $G$-orbits, under simultaneous
conjugation, in $X_j \times X_j$ of pairs $(\tau,\tau^\prime)$ that do not
satisfy (i) of Lemma \ref{invlemma}.

Take such a pair with characteristic partition 
$\{H_1,H_2,\cdots,H_t\}$.
The proofs of Lemmas \ref{lemma2} and \ref{invlemma} show that restricted to 
each $H_i$, both 
$\tau$ and $\tau^\prime$ have the same length, say, $j_i$ and further that
$|H_i| - 2j_i \leq 1$.
Thus the numerical characteristic partition of $(\tau,\tau^\prime)$ is
a partition of $N$ with $N-2j$ odd 
parts which is clearly an invariant of the orbit of $(\tau,\tau^\prime)$.
Corollary \ref{corr} implies that this invariant is complete.
In other words, two elements of $X_j \times X_j$ with the same numerical
characteristic partition are in the same $G$-orbit.

Conversely, given a partition of $N$ with $N-2j$ odd parts, take any
corresponding (set) partition, say $\{H_1,\cdots,H_t\}$, of $[N]$ and for each 
$H_i$ of cardinality
at least $2$, choose involutions $\tau_i,\tau^\prime_i \in \Sigma_{H_i}$
as in Case II or Case III of Lemma \ref{lemma2}, according as $\#H_i$ is
odd or even.
Let $\tau = \tau_1\tau_2\cdots \tau_t$ and $\tau^\prime = 
\tau_1^\prime\tau_2^\prime\cdots \tau_t^\prime$.
Then Lemma \ref{invlemma} implies that $(\tau,\tau^\prime)$  does not satisfy 
(i) and then clearly its numerical partition is the one we started with.\qed

We now have all the ingredients in place to finish the proof of the main 
theorem.
\bigskip

\noindent
{\bf Proof of Theorem \ref{main}.} It follows from Propositions \ref{disjoint} 
and \ref{comm} that $End_G(A) = \oplus_j End_G(V_j)$ and is commutative so that
any irreducible representation of $G$ occurs in $A$ at most once. Now 
Proposition \ref{dimn} implies that the dimension of $End_G(A)$ is equal to
the number of partitions of $N$ which is the number of conjugacy classes and
hence irreducible representations of $\Sigma_N$, and so each irreducible 
representation of $\Sigma_N$ does indeed feature in $A$.
\qed

\begin{Remark} This theorem may be regarded as a `structural explanation' of
the numerical coincidence arising out of the R-S-K correspondence - see 
Exercise 6 of \S 4.3 of \cite{Flt} - that the number of involutions in a 
symmetric group is equal to the
sum of the degrees of its irreducible representations.
\end{Remark}

\end{document}